\documentclass[12pt,leqno]{article}
\usepackage{amsfonts,amsthm,amsmath}
\usepackage{url}
\newtheorem{thm}{Theorem}
\newtheorem{prop}[thm]{Proposition}

\newtheorem{cor}[thm]{Corollary}
\theoremstyle{remark}
\newtheorem{rem}[thm]{Remark}

\newcommand{\ZZ}{\mathbb{Z}}
\newcommand{\cC}{\mathcal{C}}
\newcommand{\allone}{\mathbf{1}}

\begin{document}

\title{Extremal Type~II $\ZZ_4$-codes constructed from
binary doubly even self-dual codes of length $40$}

\author{
Masaaki Harada\thanks{
Research Center for Pure and Applied Mathematics,
Graduate School of Information Sciences,
Tohoku University, Sendai 980--8579, Japan.
email: mharada@m.tohoku.ac.jp.
This work was partially carried out
at Yamagata University.}}

\maketitle

\begin{abstract}
In this note, we demonstrate that
every binary doubly even self-dual code
of length $40$ can be realized as the residue code of some 
extremal Type~II $\ZZ_4$-code.
As a consequence, it is shown that
there are at least $94356$ inequivalent
extremal Type~II $\ZZ_4$-codes of length $40$.
\end{abstract}

\section{Introduction}\label{Sec:1}

Let $\ZZ_4\ (=\{0,1,2,3\})$ denote the ring of integers
modulo $4$.
A {\em $\ZZ_{4}$-code} $\cC$ of length $n$
is a $\ZZ_{4}$-submodule of $\ZZ_{4}^n$.
Two $\ZZ_4$-codes are {\em equivalent} if one can be obtained from the
other by permuting the coordinates and (if necessary) changing
the signs of certain coordinates.
A code $\cC$ is {\em self-dual} if $\cC=\cC^\perp$, where
the dual code $\cC^\perp$ is defined as
$\{ x \in \ZZ_{4}^n \mid x \cdot y = 0$ for 
all $y \in \cC\}$
under the standard inner product $x \cdot y$.
The {\em Euclidean weight} of a codeword $x=(x_1,\ldots,x_n)$ of $\cC$ is
$n_1(x)+4n_2(x)+n_3(x)$, where $n_{\alpha}(x)$ denotes
the number of components $i$ with $x_i=\alpha$ $(\alpha=1,2,3)$.
A $\ZZ_4$-code $\cC$ is {\em Type~II}
if $\cC$ is self-dual and
the Euclidean weights of all codewords of $\cC$ are divisible
by 8~\cite{Z4-BSBM} and~\cite{Z4-HSG}.
This is a remarkable class of self-dual $\ZZ_4$-codes
related to even unimodular lattices.
A Type~II $\ZZ_4$-code of length $n$ exists if and 
only if $n \equiv 0 \pmod 8$.
The {\em minimum Euclidean weight} $d_E$ of $\cC$ is the smallest
Euclidean
weight among all nonzero codewords of $\cC$.
The minimum Euclidean
weight $d_E$ of a Type~II $\ZZ_4$-code
of length $n$ is bounded by
$d_E \le 8 \lfloor n/24 \rfloor +8$~\cite{Z4-BSBM}.
A Type~II $\ZZ_4$-code meeting this bound with equality is called
{\em extremal}.

The {\em residue} code $\cC^{(1)}$ of a $\ZZ_4$-code $\cC$
is the binary code
$\cC^{(1)}= \{ c \bmod 2 \mid  c \in \cC \}$.
If $\cC$ is self-dual, then $\cC^{(1)}$ is a binary
doubly even code~\cite{Z4-CS}.
If $\cC$ is Type~II, then $\cC^{(1)}$ contains
the all-ones vector $\allone$~\cite{Z4-HSG}.
It follows that
there is a Type~II $\ZZ_4$-code $\cC$ with $\cC^{(1)}=B$
for a given binary doubly even code $B$ containing $\allone$
(see~\cite{Z4-PLF}).
However, it is not known in general whether there is an extremal
Type~II $\ZZ_4$-code $\cC$ with $\cC^{(1)}=B$ or not.

A binary doubly even self-dual code of length $n$ exists if and
only if $n \equiv 0 \pmod 8$, and
the minimum weight
$d$ of a binary doubly even self-dual code of length $n$
is bounded by $d \le 4 \lfloor n/24 \rfloor +4$.
A binary doubly even self-dual code
meeting this bound with equality is called extremal.
Two binary codes $B$ and $B'$ are equivalent
if $B$ can be obtained from $B'$ by permuting the coordinates.
The classification of binary doubly even self-dual codes
has been done for lengths up to $40$
(see~\cite{BHM40}).
For every binary doubly even self-dual code $B$ of length $24$,
there is an extremal Type~II $\ZZ_4$-code $\cC$ with 
$\cC^{(1)}=B$~\cite[Postscript]{CS97} (see also~\cite{HLM}).
In addition, 
for every binary doubly even self-dual code $B$
of length $32$, there is an extremal Type~II $\ZZ_4$-code $\cC$ 
with $\cC^{(1)}=B$~\cite{H32}.

In this note, this work is extended to length $40$.
We demonstrate that
there is an extremal Type~II $\ZZ_4$-code $\cC$ 
with $\cC^{(1)}=B$
for every binary doubly even self-dual code $B$
of length $40$. 
As a consequence, it is shown that
there are at least $94356$ inequivalent
extremal Type~II $\ZZ_4$-codes of length $40$.
In addition, our result implies that
there is an extremal Type~II $\ZZ_4$-code $\cC$ 
with $\cC^{(1)}=B$
for every binary doubly even self-dual code $B$
of length $n \in \{8,16,24,32,40\}$. 
Also, there is an extremal Type~II $\ZZ_4$-code $\cC$ 
with $\cC^{(1)}=B$
for every binary extremal doubly even self-dual code $B$
of length $n \in \{8,16,24,32,40,48\}$. 

All computer calculations in this note
were done by {\sc Magma}~\cite{Magma}.

\section{Extremal Type~II $\ZZ_4$-codes of length 40}

\subsection{Construction method}\label{SS:method}
We review the method for constructing
Type~II $\ZZ_4$-codes, which was given in~\cite{Z4-PLF}.
Let $B$ be a binary doubly even self-dual code of length $n$.
Let $I_n$ denote the identity matrix of order $n$ and
let
\[
\tilde{I_n}=
\left(\begin{array}{cccccc}
1       & \cdots   & 1 \\
0       &          &   \\
\vdots  &   I_{n-1}&   \\
0       &          &   \\
\end{array}\right).
\]
Without loss of generality, we may assume that $B$ has
generator matrix of the following form:
\begin{equation}  \label{Eq:G1}
G_1=
\left(\begin{array}{cc}
  A_1    &   \tilde{I_n} \\
\end{array}\right),
\end{equation}
where 
$A_1$ is an $n/2 \times n/2$ matrix which has the property that
the first row is  $\allone$.
Then we have a generator matrix of 
a Type~II $\ZZ_4$-code $\cC$ as follows:
\begin{equation}\label{Eq:Z4-G}
\left(\begin{array}{cccccc}
&A_1 &    & \tilde{I_n}+2A_2 \\
\end{array}\right),
\end{equation}
where $A_2$ is an $n/2 \times n/2$ $(1,0)$-matrix and
we regard the matrices as matrices over $\ZZ_4$.
Here,  we can choose freely the entries above the diagonal
elements and the $(1,1)$-entry of $A_2$, and the rest 
is completely determined from the property that $\cC$
is Type~II.
Since any Type~II $\ZZ_4$-code is equivalent to some 
Type~II $\ZZ_4$-code
containing $\allone$~\cite{Z4-HSG},
without loss of generality, we may assume that the first row
of $A_2$ is the zero vector.
This reduces 
our search space for finding extremal Type~II $\ZZ_4$-codes.
It is the aim of this work to find a $20 \times 20$ $(1,0)$-matrix
$A_2$ such that the matrix of form~\eqref{Eq:Z4-G} generates an
extremal Type~II $\ZZ_4$-code from a generator matrix
of form~\eqref{Eq:G1} 
for a given binary doubly even self-dual code of length $40$.

\subsection{Extremal Type~II $\ZZ_4$-codes of length 40}

There are $94343$ inequivalent binary doubly even self-dual 
codes of length $40$~\cite{BHM40}.
Let $B$ be one of the $94343$ binary codes.
Without loss of generality, we may assume that $B$ has
generator matrix of form~\eqref{Eq:G1}.
In the above method,
we explicitly found a $20 \times 20$ $(1,0)$-matrix
$A_2$ such that the matrix of form~\eqref{Eq:Z4-G} generates an
extremal Type~II $\ZZ_4$-code $\cC$.
Note that $\cC^{(1)}=B$.
This was done for all the $94343$ binary doubly even self-dual codes.
Hence, we have the following:

\begin{prop}
Let $B$ be a binary doubly even self-dual code of length $40$.
Then there is an extremal Type~II $\ZZ_4$-code $\cC$ with $\cC^{(1)}=B$.
\end{prop}

\begin{rem}
The extremality of the code was verified as follows.
Let $\cC$ be a Type~II $\ZZ_4$-code of length $40$.
The following lattice
\[
A_{4}(\cC) = \frac{1}{2}
\{(x_1,\ldots,x_n) \in \ZZ^n \mid
(x_1 \bmod 4,\ldots,x_n \bmod 4)\in \cC\}
\]
has minimum norm $4$ if and only if $\cC$ is extremal~\cite{Z4-BSBM}.
Instead of calculating the minimum Euclidean weight of $\cC$,
we calculated the minimum norm of $A_4(\cC)$.
This speeded up the calculations by {\sc Magma}~\cite{Magma}
considerably.
As an example, for some five extremal Type~II $\ZZ_4$-codes,
the calculations for the minimum Euclidean weights
took about $1223$ minutes, but
the calculations for the minimum norms
took about $3$ seconds only,
using a single core of a PC Intel i7 4 core processor.
\end{rem}

By the above proposition, 
$94343$ extremal Type~II $\ZZ_4$-codes are constructed
from the $94343$ inequivalent binary doubly even self-dual 
codes of length $40$.
The $94343$ extremal Type~II $\ZZ_4$-codes are
inequivalent, since their residue codes are inequivalent.
Generator matrices for the $94343$ codes
can be written in the form
$\left(\begin{array}{cc}
  I_{20}   &   M \\
\end{array}\right)$,
where $M$ can be obtained electronically from
\url{http://www.math.is.tohoku.ac.jp/~mharada/Paper/Z4-40-II.txt}.


For $m=7,8,\ldots,19$,
an extremal Type~II $\ZZ_4$-code $\cC$ of length $40$
such that the residue code $\cC^{(1)}$ has dimension $m$
is known~\cite{H32}.
Hence, we have the following:

\begin{cor}
There are at least $94356$ inequivalent
extremal Type~II $\ZZ_4$-codes of length $40$.
\end{cor}


As described above, 
for every binary doubly even self-dual code $B$ of length $24$
(resp.\ $32$),
there is an extremal Type~II $\ZZ_4$-code $\cC$ with
$\cC^{(1)}=B$~\cite[Postscript]{CS97} (resp.\ \cite{H32}).
Hence, we have the following:

\begin{cor}
Suppose that $n \in \{8,16,24,32,40\}$.
Let $B$ be a binary doubly even self-dual code of length $n$.
Then there is an extremal Type~II $\ZZ_4$-code $\cC$
with $\cC^{(1)}=B$.
\end{cor}

It is known that the binary extended quadratic residue code $QR_{48}$
of length $48$ is the unique binary extremal doubly even
self-dual code of that length.
The binary code $QR_{48}$ is the residue code of the extended lifted 
quadratic  residue $\ZZ_4$-code of length $48$, which is
an extremal Type~II $\ZZ_4$-code~\cite{Z4-BSBM}.
Hence, we have the following:

\begin{cor}
Suppose that $n \in \{8,16,24,32,40,48\}$.
Let $B$ be a binary extremal doubly even self-dual code of length $n$.
Then there is an extremal Type~II $\ZZ_4$-code $\cC$
with $\cC^{(1)}=B$.
\end{cor}

In this note, $94343$ inequivalent extremal Type~II $\ZZ_4$-codes
$\cC_i$ of length $40$ were constructed explicitly
($i=1,2,\ldots,94343$).
It is a worthwhile problem to determine whether extremal
even unimodular lattices $A_4(\cC_i)$ ($i=1,2,\ldots,94343$)
are isomorphic or not.

\bigskip
\noindent
{\bf Acknowledgment.} 
This work was supported by JSPS KAKENHI Grant Number 15H03633.


\end{document}